# Minimal $R_1$, minimal regular and minimal presober topologies

M.L. Colasante and D. Van der Zypen

**Abstract**

By means of filters, minimal $R_1$ and minimal regular topologies are characterized on suitable intervals consisting of non-trivial $R_0$ topologies.

**AMS(MOS) subject classifications.** 54A10, 54D10, 54D25

## 1    Introduction

The family $LT(X)$ of all topologies definable on a set $X$ partially ordered by inclusion is a complete, atomic lattice in which the meet of a collection of topologies is their intersection, while the join is the topology with their union as a subbase. There has been a considerable amount of interest in topologies which are minimal in this lattice with respect to certain topological properties (see for instance [1], [2], [3], [4], [5], [8], [9], [10], [11], [12], [13], [15], [18]).

Given a topological property $P$ (like a separation axiom) and given a family $\mathcal{S}$ of members of $LT(X)$, then $\tau \in \mathcal{S}$ is said to be *minimal* $P$ in $\mathcal{S}$ if $\tau$ satisfies $P$ but no member of $\mathcal{S}$ which is strictly weaker than $\tau$ satisfies $P$. It is well known that a $T_2$-topology on an infinite set $X$ is minimal $T_2$ in $LT(X)$ iff every open filter on $X$ with a unique adherent point is convergent ([3]). Also, a regular $T_1$-topology is minimal regular in $LT(X)$ iff every regular filter on $X$ with a unique adherent point is convergent ([4]). These are characterization of minimal topologies satisfying separation axioms above $T_1$, and thus topologies in the lattice $\mathcal{L}_1 = \{\tau \in LT(X) : \mathcal{C} \leq \tau \leq 2^X\}$, where $\mathcal{C}$ denotes the cofinite topology (*i.e* the minimal $T_1$-topology on $X$) and $2^X$ denotes the powerset of $X$. Some separation axioms independent of $T_1$ (even independent of $T_0$) are vacuously satisfied by the indiscrete topology, thus the study of minimal topologies in $LT(X)$ satisfying such properties becomes trivial. This is the case of the $R_1$ and regularity (not necessarily $T_1$) separation axioms. The purpose of this paper is to show that, by restricting to suitable intervals $\mathcal{L}_\rho$ of $LT(X)$, associated each to a non-trivial $R_0$-topology $\rho$, then minimal regular and minimal $R_1$ topologies in $\mathcal{L}_\rho$ can be characterized in terms of filters. For instance, we prove in section 3 that an $R_1$-topology in $\mathcal{L}_\rho$ is minimal $R_1$ iff every open filter on $X$, for which the set of adherent points coincides with a point closure, is convergent, and that a regular topology in $\mathcal{L}_\rho$ is minimal regular iff every regular filter on $X$, for which the set of adherent points coincides with a point closure, is convergent. The characterizations for minimal $T_2$ and minimal regular topologies mentioned at the beginning of this paragraph are immediate corollary of our results in case $\rho$ is a $T_1$-topology. Additionally, we consider in last section another topological property independent of $T_0$, namely the presober property, and show that there are not minimal presober topologies in $\mathcal{L}_\rho$.

## 2 Preliminaries and notations

A topology $\tau \in LT(X)$ is said to be an *Alexandroff topology* if it is closed under arbitrary intersection. J. Steprans and S. Watson [16] attributed this notion to both Alexandroff and Tucker, and called them *AT* topologies. This class of topologies is specially relevant for the study of non-$T_1$ topologies. Note that the only $T_1$ Alexandroff topology is the discrete topology. Among the characterizations known for *AT* topologies, we recall the one related with the specialization preorder: $\tau \in LT(X)$ is *AT* iff it is the finest topology on $X$ consistent with the specialization preorder, *i.e.* the finest topology giving the preorder $\leq_\tau$ satisfying $x \leq_\tau y$ iff $x$ belongs to the $\tau$-closure of $\{y\}$. This preorder characterizes the $T_0$ property (for every two points there is an open set containing one an only one of the points) in the sense that a topology $\tau$ is $T_0$ iff the preorder $\leq_\tau$ is a partial order.

By identifying a set with its characteristic function, $2^X$ can be endowed with the product topology of the Cantor cube $\{0,1\}^X$. It was proved in [17] that a topology $\tau$ on $X$ is *AT* iff it is closed when viewed as a subset of $2^X$. Moreover, it was proved there that the closure $\overline{\tau}$ of $\tau$ in $2^X$ is also a topology, and therefore it is the smallest *AT* topology containing $\tau$.

By $cl_\tau(A)$ we denote the $\tau$-closure of a set $A$. If $A = \{x\}$, we use $cl_\tau(x)$ instead of $cl_\tau(\{x\})$, and refer to it as a point closure. The $\tau$-kernel of a set $A \subseteq X$, denoted by $\ker_\tau(A)$, is the intersection of all open sets containing $A$. For any $x \in X$, we denote $\ker_\tau(\{x\}) = \ker_\tau(x)$. It is obvious that $x \in cl_\tau(y)$ iff $y \in \ker_\tau(x)$. A set $A$ is said to be $\tau$-*kernelled* (or just kernelled) if $A = \ker_\tau(A)$. Equivalently, $A$ is kernelled iff $A = \bigcup_{x \in A} \ker_\tau(x)$. The family of all kernelled subsets of $X$ is closed under arbitrary unions and intersections, so it is an *AT* topology. Moreover, it coincides with $\overline{\tau}$. In fact, since every open set is kernelled and $\overline{\tau}$ is the smallest $AT$ topology containing $\tau$, then every member of $\overline{\tau}$ is kernelled. On the other hand, since $\overline{\tau}$ is closed under arbitrary intersections and it contains $\tau$, then every kernelled set belongs to $\overline{\tau}$. Thus, $\overline{\tau}$ is the topology on $X$ generated by the family $\{\ker_\tau(x) : x \in X\}$. In particular, $A \subseteq X$ is $\overline{\tau}$-closed iff $A = \bigcup_{x \in A} cl_\tau(x)$. Note that, since $\tau$ is $T_1$ iff every subset of $X$ is kernelled, then $\tau$ is $T_1$ iff $\overline{\tau} = 2^X$.

In what follows $\mathcal{N}_\tau(x)$ denotes the filter base of $\tau$-neighborhoods of $x \in X$. A filter $\mathcal{F}$ on $X$ is said to be $\tau$-convergent to a point $x \in X$ if $\mathcal{F} \supseteq \mathcal{N}_\tau(x)$. By $adh_\tau \mathcal{F}$ we denote the set of adherent points of $\mathcal{F}$ (i.e. $adh_\tau \mathcal{F} = \bigcap_{F \in \mathcal{F}} cl_\tau(F)$). Since $adh_\tau \mathcal{F}$ is a closed set, then it contains the $\tau$-closure of all its points. It is immediate that if $\mathcal{F}$ is $\tau$-convergent to $x$, then $\mathcal{F}$ is $\tau$-convergent to every $y \in cl_\tau(x)$. A filter $\mathcal{F}$ is said to be $\tau$-open if $F \in \tau$ for all $F \in \mathcal{F}$, and $\mathcal{F}$ is said to be $\tau$-regular if it is $\tau$-open and for every $F \in \mathcal{F}$ there exists $F' \in \mathcal{F}$ such that $cl_\tau(F') \subseteq F$. Thus, a $\tau$-regular filter is equivalent to a $\tau$-closed filter. A filter on $X$ is said to be ultrafilter if it is a maximal filter.

For definitions and notations not given here, we refer the reader to [19].

## 3    Minimal $R_1$ and minimal regular topologies in $\mathcal{L}_\rho$

In this section, we restrict our attention to suitable intervals consisting of $R_0$ topologies, and give characterizations of minimal $R_1$ and minimal regular topologies on those intervals. Recall that a topology $\tau \in LT(X)$ is said to be:

($R_0$) *if for all* $x, y \in X$, $x \in cl_\tau(y)$ *iff* $y \in cl_\tau(x)$, *thus* $\tau$ *is* $R_0$ *iff the point closures form a partition of* $X$. [14]

($R_1$) *if for all* $x, y \in X$ *with* $cl_\tau(x) \neq cl_\tau(y)$, *there are disjoint open sets* separating $cl_\tau(x)$ and $cl_\tau(y)$. [7]

($Regular$) *if for each* $V \in \tau$ *and each* $x \in V$ *there exists* $U \in \tau$ *such that* $x \in U \subseteq cl_\tau(U) \subseteq V$.

The separation axioms $R_0$ and $R_1$ are also denoted as $S_1$ and $S_2$, respectively ([6]). We use in this paper the most common notations $R_0$ and $R_1$. It is easy to show that $Regularity \Rightarrow R_1 \Rightarrow R_0$, and that none of the implications can be reversed. Moreover, $\tau$ is $T_1$ iff $\tau$ is $R_0$ and $T_0$, and $\tau$ is $T_2$ iff $\tau$ is $R_1$ and $T_0$.

Examples of topologies which are regular non-$T_0$ (thus, regular non-$T_1$) abound. For instance, if $\mathcal{P}$ denotes any non-trivial partition of a set $X$, then the associated partition topology $\tau_\mathcal{P}$, defined as the topology having as open sets the unions of elements of $\mathcal{P}$ together with the empty set, is a regular topology which is not $T_0$. On the other hand, if a topological space satisfies any of the $(R)$ properties stated above and one doubles the space by taking the product of $X$ with the two point indiscrete space, then the resulting space is not longer $T_0$ but it satisfies the same $(R)$ properties as did the original space.

The following characterizations, which are straightforward to prove, are used throughout the paper without explicitly mentioning them.

**Lemma 3.1** *Let $\tau \in LT(X)$. Then*

*(i) $\tau$ is $R_0$ iff $cl_\tau(x) = \ker_\tau(x)$ for all $x \in X$, iff $cl_\tau(x) \subseteq V$, for all $V \in \tau$ and $x \in V$.*

*(ii) $\tau$ is $R_1$ iff $\tau$ is $R_0$ and for all $x, y \in X$ such that $y \notin cl_\tau(x)$ there are disjoint open sets separating $x$ and $y$.*

*(iii) $\tau$ is $R_1$ iff $\tau$ is $R_0$ and $adh_\tau \mathcal{N}_\tau(x) = cl(x)$, for all $x \in X$.*

To each $\rho \in LT(X)$ we associate the interval

$$\mathcal{L}_\rho = \{\tau \in LT(X) : at(\rho) \leq \tau \leq \overline{\rho}\}$$

where $at(\rho)$ denotes the topology on $X$ generated by the sets $\{X \setminus cl_\rho(H) : H$ is a finite subset of $X\}$, and $\overline{\rho}$ is the closure of $\rho$ in $2^X$.

Note that, if $\rho$ is any $T_1$-topology, then $at(\rho) = \mathcal{C}$ and $\overline{\rho} = 2^X$. In this case, $\mathcal{L}_\rho$ is precisely the lattice $\mathcal{L}_1$ of all $T_1$ topologies on $X$.

**Lemma 3.2** *Let $\rho \in LT(X)$. Then $cl_{at(\rho)}(x) = cl_\rho(x) = cl_{\overline{\rho}}(x)$, for every $x \in X$.*

**Proof.** Let $x \in X$. Since a set is $\overline{\rho}$-closed iff it is union of $\rho$-closed sets, then $cl_\rho(x) \subseteq cl_{\overline{\rho}}(x)$. On the other hand, $cl_\rho(x)$ is an $at(\rho)$-closed set, and thus $cl_{at(\rho)}(x) \subseteq cl_\rho(x)$. Since $at(\rho) \subseteq \rho \subseteq \overline{\rho}$, then $cl_{\overline{\rho}}(x) \subseteq cl_\rho(x) \subseteq cl_{at(\rho)}(x)$. From this we have the result. ∎

**Corollary 3.3** *Let $\tau, \rho \in LT(X)$. Then $\tau \in \mathcal{L}_\rho$ iff $cl_\rho(x) = cl_\tau(x)$, for every $x \in X$.*

**Proof.** If $\tau \in \mathcal{L}_\rho$ and $x \in X$, then Lemma 3.2 implies that $cl_\tau(x) = cl_\rho(x)$. Conversely, suppose $cl_\rho(x) = cl_\tau(x)$, for every $x \in X$. It is immediate that $at(\rho) \leq \tau$. Note that $\ker_\rho(x) = \ker_\tau(x)$, thus if $V \in \tau$ then $V = \bigcup_{x \in V} \ker_\rho(x) = \bigcup_{x \in V} \ker_\tau(x)$ is a $\overline{\rho}$-open set. Therefore $at(\rho) \leq \tau \leq \overline{\rho}$. ∎

Corollary 3.3 can be stated as follows: $\tau \in \mathcal{L}_\rho$ iff $\tau$ has the same preorder of specialization as $\rho$. Thus, when one refers to the $\tau$-closure of $x \in X$, for any $\tau \in \mathcal{L}_\rho$, there is not need for specifying the topology. We will often write $cl(x)$ without further comment. It is clear that the topologies on $\mathcal{L}_\rho$ share the topological properties defined in terms of point closures. In particular $\tau \in \mathcal{L}_\rho$ is $R_0$ iff $\rho$ is $R_0$. Note that the property $R_1$ is expansive in $\mathcal{L}_\rho$ (i.e. if $\tau \in \mathcal{L}_\rho$ is $R_1$, then $\tau'$ is $R_1$ for all $\tau' \in \mathcal{L}_\rho$ finer than $\tau$).

In ([6]) it was proved that the properties $R_0$, $R_1$ and regularity coincide for $AT$ topologies. Thus, $\overline{\rho}$ is $R_0$ iff $\overline{\rho}$ is $R_1$ iff $\overline{\rho}$ is regular. If we start with an $R_0$-topology $\rho$ on $X$, it is immediate that there exists at least a regular topology (so at least an $R_1$-topology) in $\mathcal{L}_\rho$. Our goal is to characterize minimal $R_1$ and minimal regular topologies in $\mathcal{L}_\rho$. Note that, if $\rho$ is $R_0$ and $X$ can be written as a finite union of disjoint point closures then, for each $x \in X$, the set $cl(x)$ is the complement of finite union of point closures, thus $cl(x) \in at(\rho)$. It follows that $at(\rho) = \rho = \overline{\rho}$ and therefore $\mathcal{L}_\rho = \{\rho\}$. To avoid triviality, from now on we assume that $\rho \in LT(X)$ is an $R_0$-topology such that $X$ can be written as infinite union of disjoint point closures (in particular, this is the case for any $T_1$-topology on an infinite set). It is worth notice that $at(\rho)$ can not be $R_1$, thus it can not be regular, since any pair of non-empty $at(\rho)$-open sets intersect. We give an example of an $R_0$ (not $T_0$) topology satisfying the above conditions.

**Example 3.4** *Let $X$ be the set of all positive integers $N$, and let $\rho$ be the topology generated by the subbase $\{\emptyset, N\backslash\{1\}, N\backslash\{2n, 2n+1\}, n \geq 1\}$. It is easy to see that $\rho$ is an $R_0$-topology which is not $T_0$, and that $N$ can be written as the infinite disjoint union of the odd integers point closure. Note that $at(\rho) = \rho$, and $\overline{\rho}$ is the topology generated by the sets $\{1\}, \{2n, 2n+1\}, n \geq 1$.*

For $x \in X$, let $\mathcal{E}(x)$ denote the family of all the subsets of $X$ not containing $x$. If $\mathcal{F}$ is any filter on $X$, then $\mathcal{E}(x) \cup \mathcal{F}$ is a topology on $X$. Given $\tau \in LT(X)$, we consider the topology $\beta = \tau \cap (\mathcal{E}(x) \cup \mathcal{F})$. Note that $\beta \leq \tau$, and $\beta = \tau$ iff $\mathcal{F} = \mathcal{N}_\tau(x)$.

Now, if $\rho$ is $R_0$ and $\tau \in \mathcal{L}_\rho$, a local base for the topology $\beta$ can be describe as follows:

$\mathcal{N}_\beta(y) = \mathcal{N}_\tau(y) \cap \mathcal{E}(x)$, for every $y \notin cl(x)$;

$\mathcal{N}_\beta(y) = \mathcal{N}_\tau(x) \cap \mathcal{F}$, for every $y \in cl(x)$.

A set $A \subseteq X$ is $\beta$-closed iff $A$ is $\tau$-closed and either $x \in A$ or $X\backslash A \in \mathcal{F}$. Thus, $cl_\tau(A) \subseteq cl_\beta(A) \subseteq cl_\tau(A) \cup cl(x)$ for all $A \subseteq X$. In particular $cl_\beta(x) = cl(x)$.

**Lemma 3.5** *Let $\tau \in \mathcal{L}_\rho$. Given $x \in X$ and a filter $\mathcal{F}$ on $X$, let $\beta = \tau \cap (\mathcal{E}(x) \cup \mathcal{F})$. Then*

*(i) $\beta$ is $R_0$ iff $\mathcal{F} \supseteq \mathcal{N}_{at(\rho)}(x)$ iff $\beta \in \mathcal{L}_\rho$.*

*(ii)* If $adh_\tau \mathcal{F} = cl(x)$, then $\beta \in \mathcal{L}_\rho$.

**Proof.** (i) It is immediate that $\mathcal{F} \supseteq \mathcal{N}_{at(\rho)}(x)$ iff $\beta \in \mathcal{L}_\rho$, and that if $\beta \in \mathcal{L}_\rho$ then $\beta$ is $R_0$. On the other hand, if $\beta$ is $R_0$ and $y \notin cl(x) = cl_\beta(x)$, then $x \notin cl_\beta(y)$. Thus $X \setminus cl_\beta(y) \in \mathcal{F}$, and this implies that $X \setminus cl(y) \in \mathcal{F}$. Since this holds for every $y \notin cl(x)$, it follows that $\mathcal{F} \supseteq \mathcal{N}_{at(\rho)}(x)$.

(ii) If $adh_\tau \mathcal{F} = cl(x)$ and $y \notin cl(x)$ then $y \notin adh_\tau \mathcal{F}$, and thus there exist $F \in \mathcal{F}$ and $V \in \mathcal{N}_\tau(y)$ such that $V \cap F = \emptyset$. Since $cl(y) \subseteq V$, then $F \subseteq X \setminus cl(y)$ and thus $X \setminus cl(y) \in \mathcal{F}$. Hence $\mathcal{F} \supseteq \mathcal{N}_{at(\rho)}(x)$. ∎

**Proposition 3.6** *Let $\tau \in \mathcal{L}_\rho$ be $R_1$. Given $x \in X$ and a filter $\mathcal{F}$ on $X$, then the topology $\beta = \tau \cap (\mathcal{E}(x) \cup \mathcal{F})$ is $R_1$ iff there exists a $\tau$-open filter $\mathcal{F}_0 \subseteq \mathcal{F}$ such that $adh_\tau \mathcal{F}_0 = cl(x)$.*

**Proof.** ($\Rightarrow$) If $\beta = \tau \cap (\mathcal{E}(x) \cup \mathcal{F})$ is $R_1$, then $adh_\beta \mathcal{N}_\beta(x) = cl(x)$. By Lemma 3.5(i), $\beta \in \mathcal{L}_\rho$. Now, since $\beta \leq \tau$, then $cl(x) \subseteq adh_\tau \mathcal{N}_\beta(x) \subseteq adh_\beta \mathcal{N}_\beta(x) = cl(x)$. Let $\mathcal{F}_0 = \mathcal{N}_\beta(x) = \mathcal{N}_\tau(x) \cap \mathcal{F}$. It is clear that $\mathcal{F}_0$ is a $\tau$-open filter contained in $\mathcal{F}$ such that $adh_\tau \mathcal{F}_0 = cl(x)$.

($\Rightarrow$) Suppose there exists a $\tau$-open filter $\mathcal{F}_0 \subseteq \mathcal{F}$ such that $adh_\tau \mathcal{F}_0 = cl(x)$. By Lemma 3.5(i), $\beta \in \mathcal{L}_\rho$. To prove that $\beta = \tau \cap (\mathcal{E}(x) \cup \mathcal{F})$ is $R_1$, let $y, z \in X$ such that $y \notin cl(z)$. We will show that $y$ and $z$ can be separated by $\beta$-open sets. Since $\tau$ is $R_1$, there exist $W_y \in \mathcal{N}_\tau(y)$ and $W_z \in \mathcal{N}_\tau(z)$ such that $W_y \cap W_z = \emptyset$. We consider two possible cases.

*Case (i).* If $x \notin cl(y)$ and $x \notin cl(z)$, then $y, z \notin cl(x)$. Choose $V_y \in \mathcal{N}_\tau(y)$ and $V_z \in \mathcal{N}_\tau(z)$ such that $x \notin V_y$ and $x \notin V_z$. Let $O_y = W_y \cap V_y$ and $O_z = W_z \cap V_z$. Then $O_y, O_z \in \tau \cap \mathcal{E}(x) \leq \beta$ and $O_y \cap O_z = \emptyset$.

*Case (ii).* If $x \in cl(y)$, then $cl(y) = cl(x) = adh_\tau \mathcal{F}_0$. Since $z \notin cl(y)$, there exists $U \in \mathcal{N}_\tau(z)$ and $F \in \mathcal{F}_0$ such that $U \cap F = \emptyset$. Take $O_y = W_y \cup F$ and $O_z = W_z \cap U$. Then it is immediate that $O_y \in \tau \cap \mathcal{F}$ and $O_z \in \tau \cap \mathcal{E}(x)$. Thus $O_y$ and $O_z$ are disjoint $\beta$-neighborhoods of $y$ and $z$, respectively. ∎

**Remark 3.7** *For any $x \in X$, the open filter $\mathcal{F} = \mathcal{N}_{at(\rho)}(x)$ satisfies $adh_{\overline{\rho}} \mathcal{F} = cl(x)$. In fact, $cl(y) = \ker(y) \in \mathcal{N}_{\overline{\rho}}(y)$ for each $y \in X$. Then, $y \in cl(x)$ implies that $x \in V$, for all $V \in \mathcal{N}_{\overline{\rho}}(y)$, and thus $y \in adh_{\overline{\rho}} \mathcal{N}_{at(\rho)}(x)$. On the other hand, if $y \notin cl(x)$ then the disjoint sets $cl(y) \in \mathcal{N}_{\overline{\rho}}(y)$ and $X \setminus cl(y) \in \mathcal{N}_{at(\rho)}(x)$ witness that $y \notin adh_{\overline{\rho}} \mathcal{N}_{at(\rho)}(x)$. Since $\overline{\rho}$ is $R_1$, above proposition implies that $\beta = \overline{\rho} \cap (\mathcal{E}(x) \cup \mathcal{N}_{at(\rho)}(x))$ is $R_1$ and hence $\beta \in \mathcal{L}_\rho$. Moreover, $\beta$ is strictly weaker than $\overline{\rho}$ since $cl(x) \in \overline{\rho}$ but $cl(x) \notin \mathcal{N}_{at(\rho)}(x)$. Therefore $\overline{\rho}$ is not the minimal $R_1$ topology in $\mathcal{L}_\rho$.*

We are now ready to prove a characterization of minimal $R_1$ in $\mathcal{L}_\rho$.

**Theorem 3.8** *Let $\tau \in \mathcal{L}_\rho$ be $R_1$. Then $\tau$ is minimal $R_1$ iff given any open filter $\mathcal{F}$ on $X$ such that $adh_\tau \mathcal{F} = cl(x)$ for some $x \in X$, then $\mathcal{F}$ is convergent (necessarily to every point of $cl(x)$).*

**Proof.** Suppose $\tau$ is minimal $R_1$ and let $\mathcal{F}$ be an open filter on $X$ such that $adh_\tau \mathcal{F} = cl(x)$, for some $x \in X$. Let $\beta = \tau \cap (\mathcal{E}(x) \cup \mathcal{F})$. By Lemma 3.5(i), $\beta \in \mathcal{L}_\rho$ and, by Proposition 3.6, $\beta$ is $R_1$. Since $\tau$ is minimal $R_1$ in $\mathcal{L}_\rho$, it must be that $\beta = \tau$, and thus $\mathcal{F} \supseteq \mathcal{N}_\tau(x)$.

Conversely, suppose every open filter $\mathcal{F}$ on $X$ such that $adh_\tau \mathcal{F} = cl(x)$, for some $x \in X$, is $\tau$-convergent and let $\tau' \in \mathcal{L}_\rho$ be an $R_1$-topology such that $\tau' \leq \tau$. Let $V \in \tau$ and $x \in V$. Since $adh_{\tau'} \mathcal{N}_{\tau'}(x) = cl(x)$, the hypothesis implies that the $\tau$-open filter $\mathcal{N}_{\tau'}(x)$ is $\tau$-convergent to $x$. Thus $\mathcal{N}_{\tau'}(x) \supseteq \mathcal{N}_\tau(x)$, and hence $V \in \mathcal{N}_{\tau'}(x)$. Since this happens for all $x \in V$, then $V \in \tau'$. Therefore $\tau = \tau'$, and this implies that $\tau$ is minimal $R_1$. ∎

Since $\tau$ is minimal $T_2$ iff $\tau \in \mathcal{L}_1$ and is minimal $R_1$, then Theorem 3.8 applied to any $T_1$-topology $\rho$ yields the following well known result on minimal $T_2$.

**Corollary 3.9** *Let $X$ be an infinite set, and let $\tau \in LT(X)$ be $T_2$. Then $\tau$ is minimal $T_2$ iff every open filter on $X$ with a unique adherent point is convergent (to that point).*

Recall that $\tau \in LT(X)$ is said to be compact if every open cover of $X$ has a finite subcover. Equivalently, $\tau$ is compact iff every filter on $X$ has an adherent point iff every ultrafilter on $X$ converges. ([19]). It is known that if $\tau$ is minimal $T_2$ then $\tau$ is regular iff it is compact ([19]). We will show that this last equivalence holds for minimal $R_1$ topologies in $\mathcal{L}_\rho$. The results given in the following lemma are well known. For the sake of completeness, we include the proofs.

**Lemma 3.10** *Let $\tau \in LT(X)$.*

*(i) If $\tau$ is $R_1$ and compact, then $\tau$ is regular.*

*(ii) If $\tau$ is regular and every open filter on $X$ has an adherent point, then $\tau$ is compact.*

**Proof.** (i) Let $\tau$ be $R_1$ and compact, and let $x \in V \in \tau$. Then for each $y \in X \setminus V$, there exist $U^y \in \mathcal{N}_\tau(x)$ and $V_y \in \mathcal{N}_\tau(y)$ such that $U^y \cap V_y = \emptyset$. Now, the family $\{V_y\}_{y \in X \setminus V}$ is an open cover of $X \setminus V$, a closed set and hence a compact set. Thus, $X \setminus V \subseteq \bigcup_{i=1}^n V_{y_i}$ for some finite collection $\{y_1, ..., y_n\}$ of points in $X \setminus V$. Let $U = \bigcap_{i=1}^n U^{y_i}$. It is immediate that $U \in \mathcal{N}_\tau(x)$ and $cl(U) \subseteq V$, which shows that $\tau$ is regular.

(ii) Let $\tau$ be regular and such that every open filter on $X$ has an adherent point. Given an ultrafilter $\mathcal{R}$ on $X$, consider the open filter $\mathcal{F} = \mathcal{R} \cap \tau$. Then $\mathcal{F}$ has an adherent point $x \in X$. Now, if $\mathcal{R}$ does not converge to $x$, there exists $V \in \mathcal{N}_\tau(x)$ such that $V \notin \mathcal{R}$ and hence $X \setminus V \in \mathcal{R}$, since $\mathcal{R}$ is ultrafilter. By regularity of $\tau$, one can choose $U \in \mathcal{N}_\tau(x)$ with $cl_\tau(U) \subseteq V$. Then $X \backslash cl_\tau(U) \supseteq X \backslash V$ and thus $X \backslash cl_\tau(U) \in \mathcal{R} \cap \tau = \mathcal{F}$. But, since $x \in adh_\tau \mathcal{F}$, it must be that $U \cap X \backslash cl_\tau(U) \neq \emptyset$, a contradiction. Thus $\mathcal{R}$ converges to $x$, and therefore $\tau$ is compact. ∎

**Proposition 3.11** *Let $\tau \in \mathcal{L}_\rho$. If $\tau$ is minimal $R_1$, then every open filter on $X$ has an adherent point.*

**Proof.** Suppose there is an open filter $\mathcal{F}$ on $X$ such that $adh_\tau \mathcal{F} = \emptyset$. For each $x \in X$ there exist $V \in \mathcal{N}_\tau(x)$ and $F \in \mathcal{F}$ such that $V \cap F = \emptyset$. In particular $V \notin \mathcal{F}$. On the other hand, since $cl(x) \subseteq V$ then $X \backslash cl(x) \supseteq X \backslash V \supseteq F$, and thus $X \backslash cl(x) \in \mathcal{F}$. This shows that $\mathcal{F} \supseteq \mathcal{N}_{at(\rho)}(x)$ and $\mathcal{F} \not\supseteq \mathcal{N}_\tau(x)$, for each $x \in X$. Now, fix $x \in X$ and let $\beta = \tau \cap (\mathcal{E}(x) \cup \mathcal{F})$. Then $\beta$ is a topology in $\mathcal{L}_\rho$ which is strictly weaker than $\tau$. We will prove that $\beta$ is $R_1$ and thus $\tau$ is not minimal $R_1$.

By Proposition 3.6, it is enough to show that $\mathcal{F}$ contains an open filter $\mathcal{F}_0$ such that $adh_\tau \mathcal{F}_0 = cl(x)$. Let $\mathcal{F}_0 = \{F \in \mathcal{F} : F \cap V \neq \emptyset,$ for all $V \in \mathcal{N}_\tau(x)\}$. It is clear that $\mathcal{F}_0$ is an open non-empty proper sub-filter of $\mathcal{F}$ and that $cl(x) \subseteq adh_\tau \mathcal{F}_0$. Now, let $y \notin cl(x)$. Since $\tau$ is $R_1$, there exist $V \in \mathcal{N}_\tau(x)$ and $W \in \mathcal{N}_\tau(y)$ such that $V \cap W = \emptyset$. On the other hand, since $y \notin adh_\tau \mathcal{F}$, there exist $U \in \mathcal{N}_\tau(y)$ and $F \in \mathcal{F}$ such that $U \cap F = \emptyset$. If $O = W \cap U$ and $G = V \cup F$, then $O \in \mathcal{N}_\tau(y)$, $G \in \mathcal{F}_0$ and $O \cap G = \emptyset$. Thus $y \notin adh_\tau \mathcal{F}_0$ and therefore $adh_\tau \mathcal{F}_0 = cl(x)$. ∎

Next result follows immediately from Lemma 3.10 and Proposition 3.11.

**Theorem 3.12** *Let $\tau \in \mathcal{L}_\rho$ be minimal $R_1$. Then $\tau$ is compact iff it is regular.*

We end this section with the characterization of minimal regular topologies in $\mathcal{L}_\rho$ announced over the introduction of this paper.

**Theorem 3.13** *Let $\tau \in \mathcal{L}_\rho$ be regular. Then $\tau$ is minimal regular iff every regular filter $\mathcal{F}$ on $X$ such that $adh_\tau \mathcal{F} = cl(x)$, for some $x \in X$, is convergent (necessarily to every point of $cl(x)$).*

**Proof.** ($\Rightarrow$) Let $\mathcal{F}$ be a $\tau$-regular filter on $X$ such that $adh_\tau \mathcal{F} = cl(x)$, for some $x \in X$, and suppose $\mathcal{F}$ does not converge. Then, there exists $U \in \mathcal{N}_\tau(x)$ such that $U \notin \mathcal{F}$, and hence $\beta = \tau \cap ((\mathcal{E}(x) \cup \mathcal{F}) \in \mathcal{L}_\rho$ is strictly weaker than $\tau$. Note that $x \in F$ for all $F \in \mathcal{F}$. Otherwise $x \notin cl_\tau(F')$ for some $F' \in \mathcal{F}$ and hence $x \notin adh_\tau \mathcal{F}$, which contradicts the hypothesis that $adh_\tau \mathcal{F} = cl(x)$. We prove that $\beta$ is regular and therefore $\tau$ is not minimal regular.

Let $V \in \beta$ and $y \in V$. If $y \in cl(x)$, then $V \in \mathcal{N}_\tau(x) \cap \mathcal{F} = \mathcal{F}$, a regular filter, and thus there exists $U \in \mathcal{F}$ such that $cl_\tau(U) \subseteq .V$. Since $x \in U$, then $cl_\beta(U) = cl_\tau(U) \subseteq V$. Now, if $y \notin cl(x) = adh_\tau \mathcal{F}$, there exist $U' \in \mathcal{N}_\tau(y)$ and $F \in \mathcal{F}$ such that $U' \cap F = \emptyset$. Choose $U \in \mathcal{N}_\tau(y)$ such that $cl_\tau(U) \subseteq V$ (this is possible since $\tau$ is regular). If $W = U \cap U'$, then $cl_\tau(W) \cap F = \emptyset$ and thus $X \backslash cl_\tau(W) \in \mathcal{F}$. It follows that $cl_\beta(W) = cl_\tau(W) \subseteq cl_\tau(U) \subseteq V$.

($\Leftarrow$) Suppose that every $\tau$-regular filter on $X$ for which the set of adherent points coincides with a point closure, is $\tau$-convergent. Let $\tau' \in \mathcal{L}_\rho$ be a regular topology such that $\tau' \leq \tau$. Fix $V \in \tau$ and $x \in V$. It is clear that $cl(x) = adh_\tau \mathcal{N}_\tau(x) \subseteq adh_\tau \mathcal{N}_{\tau'}(x) \subseteq adh_{\tau'} \mathcal{N}_{\tau'}(x) = cl(x)$. Since $\mathcal{N}_{\tau'}(x)$ is a $\tau'$-regular filter, then $\mathcal{N}_{\tau'}(x)$ is a $\tau$-regular filter. By hypothesis, $\mathcal{N}_{\tau'}(x)$ is $\tau$-convergente, i.e. $\mathcal{N}_\tau(x) \subseteq \mathcal{N}_{\tau'}(x)$. Since this holds for every $x \in V$, then $V \in \tau'$ and thus $\tau' = \tau$. Therefore $\tau$ is minimal regular in $\mathcal{L}_\rho$. ∎

**Corollary 3.14** *A regular and $T_1$-topology on $X$ is minimal regular iff every regular filter on $X$ with a unique adherent point is convergent.*

**Proof.** Apply Theorem 3.13 to any $T_1$-topology $\rho$. ∎

# 4 Presober topologies in $\mathcal{L}_\rho$

In this last section we consider a topological property known as presoberty, which is strictly weaker than $R_1$, and show that there are not minimal presober topologies in $\mathcal{L}_\rho$. As in previous section, we assume $\rho \in LT(X)$ is any $R_0$-topology such that $X$ can be written as infinite union of disjoint point closures.

**Definition 4.1** *A non-empty closed subset $C$ of $X$ is said reducible if there are non-empty, proper closed subsets $C_1, C_2$ of $C$, such that $C = C_1 \cup C_2$. Otherwise $C$ is irreducible. By convention, $\emptyset$ is neither reducible nor irreducible.*

Every point closure is irreducible. If $C$ is an irreducible closed set then it may be the case that it is the point closure of some point $x$. If so, $x$ is called a generic point of $C$

**Definition 4.2** *A topology is said to be presober iff each irreducible closed set has at least one generic point.*

In case that every irreducible closed subset of a space has a unique generic point, the topology is said to be *sober*. Soberty is thus a combination of two properties: the existence of generic points and their uniqueness. It is straightforward to see that the generic points in a topological space are unique iff the space satisfies the $T_0$ separation axiom. Thus, a topology is sober precisely when it is $T_0$ and presober.

In any $T_2$-topology, the irreducible closed sets are the singleton, so $T_2$ implies soberty. The cofinite topology on an infinite set is an example of a $T_1$-topology which is not sober, so it is also an example of a $T_0$ and not presober topology.

**Proposition 4.3** *Every $R_1$-topology $\tau \in LT(X)$ is presober.*

**Proof.** Let $\tau \in LT(X)$ be $R_1$, and let $C \subseteq X$ be closed. Let $x, y \in C$ with $x \neq y$. Then $cl_\tau(x)$, $cl_\tau(y) \subseteq C$. If $y \notin cl_\tau(x)$, there exist disjoint open sets $U \in \mathcal{N}_\tau(x)$ and $V \in \mathcal{N}_\tau(y)$ such that $cl_\tau(x) \subseteq U$ and $cl_\tau(y) \subseteq V$. Let $C_1 = C \cap X \backslash U$ and $C_2 = C \cap X \backslash V$. Then $C_1, C_2$ are non-empty proper closed subsets of $C$ such that $C_1 \cup C_2 = C$, and thus $C$ is reducible. It follows that an irreducible closed set must be a point closure, and hence $\tau$ is presober. ∎

Presoberty does not imply $R_1$, as the next example shows.

**Example 4.4** *Let $X$ be a set with cardinality $\geq 3$, and let $a, b \in X$ with $a \neq b$. Let $\tau$ be the topology $\{G \subseteq X : \{a,b\} \subseteq G\} \cup \{\emptyset\}$. So, a set $C$ is closed iff $C \cap \{a.b\} = \emptyset$ or $C = X$. It is clear that every $x \notin \{a,b\}$ is closed. If $C$ is non-empty, closed proper subset of $X$, then $C$ is irreducible iff it is a singleton $x \notin \{a,b\}$, since otherwise $C = \{x\} \cup (C - \{x\})$ for any $x \in C$, and both $\{x\}$ and $C - \{x\}$ are closed and non-empty. Also $X$ is itself irreducible since it is a point closure, $X = cl_\tau(a) = cl_\tau(b)$. Thus the irreducible closed sets are the point closures, and so $\tau$ is presober. But $\tau$ is not $R_1$ since given any $x \notin \{a,b\}$, then $cl_t(x)$ and $cl_\tau(a)$ can not be separated by disjoint open sets. Note that $\tau$ is an Alexandroff not $T_0$ topology on $X$.*

**Proposition 4.5** *The presober property is expansive in $\mathcal{L}_\rho$ (i.e. if $\tau \in \mathcal{L}_\rho$ is presober, then $\tau'$ is presober for all $\tau' \in \mathcal{L}_\rho$ finer than $\tau$).*

**Proof.** Let $\tau \in \mathcal{L}_\rho$ be presober, and let $\tau' \in \mathcal{L}_\rho$ with $\tau' \leq \tau$. Given a not empty $\tau'$-closed subset $A$ of $X$, let $B = cl_\tau(A)$. If $B$ is $\tau$-reducible and $F$ and $G$ are two not empty, $\tau$-closed, proper subsets of $B$ such that $B = F \cup G$, then $F_1 = (A \cap F)$ and $G_1 = (A \cap G)$ are two not empty, $\tau'$-closed, proper subsets of $A$ such that $A = F_1 \cup G_1$. Hence $A$ is $\tau'$-reducible. On the other hand, if $B$ is $\tau$-irreducible there exists $b \in B$ such that $B = cl(b)$, since $\tau$ is presober. Note that $b \in A$. Otherwise, if $a$ is any point of $A$, then $b \notin cl(a)$. Since $\tau$ is $R_0$, one has that $cl(a) \cap cl(b) = \emptyset$ and hence $cl(a) \cap B = \emptyset$, which contradicts the fact that $A \subseteq B$. Thus $cl(b) = B \subseteq A$, and it follows that $A$ is $\tau$-closed. We have proved that the $\tau'$-irreducible subsets of $X$ are $\tau$-irreducible, and thus a point closure. Therefore $\tau'$ is presober. ∎

Since $\rho$ is $R_0$, then $\overline{\rho}$ is $R_1$ and thus it is presober. Therefore, there exists at least a presober member of $\mathcal{L}_\rho$. On the other hand, $at(\rho)$ is not presober since a proper subset of $X$ is $at(\rho)$-closed iff it is finite union of disjoint point closure sets, and thus $X$ is $at(\rho)$-irreducible, but $X$ is not a point closure. Thus, $at(\rho)$ is an example of an $R_0$-topology which is not presober. We will prove that there are not minimal presober topologies in $\mathcal{L}_\rho$.

Given $\tau \in \mathcal{L}_\rho$, $x \in X$ and $\mathcal{F}$ a filter on $X$, lets consider the topology $\beta = \tau \cap (\mathcal{E}(x) \cup \mathcal{F})$.

**Lemma 4.6** *Let $\tau \in \mathcal{L}_\rho$ be presober, and let $A \subseteq X$ be $\beta$-closed. If $A$ is $\tau$-reducible then it is also $\beta$-reducible.*

**Proof.** Let $A \subseteq X$ be $\beta$-closed and $\tau$-reducible, and let $F$ and $G$ be non-empty $\tau$-closed proper subsets of $A$ such that $A = F \cup G$. Then either $x \in A$ or $X \backslash A \in \mathcal{F}$. If $X \backslash A \in \mathcal{F}$ or $x \in F \cap G$ then $F$ and $G$ are $\beta$-closed and therefore $A$ is $\beta$-reducible. Thus, we just need to consider the case when $x$ belongs to only one of the sets $F$ or $G$.

Suppose $x \in F \backslash G$ (the case $x \in G \backslash F$ is similar). Then, it is clear that $F$ is $\beta$-closed. Moreover, since $x \notin G$ and since $\tau$ is $R_0$, it must be that $cl(x) \cap G = \emptyset$ (if $y \in cl(x) \cap G$ then $x \in cl(y) \subseteq G$). Write $A = F \cup \{cl(x) \cup G\}$. If $F \backslash \{cl(x) \cup G\} \neq \emptyset$, then $F$ and $cl(x) \cup G$ are non-empty $\beta$-closed proper subsets of $A$, and thus $A$ is $\beta$-reducible. If $F \backslash \{cl(x) \cup G\} = \emptyset$, we distinguish the following cases:

(i) $G$ is $\tau$-irreducible. In this case $G = cl(g)$ for some $g \in G$, since $\tau$ is presober. Thus $A = cl(x) \cup cl(g)$ and therefore $A$ is $\beta$-reducible.

(ii) $G$ is $\tau$-reducible. Then there exist $G_1$ and $G_2$, non empty $\tau$-closed proper subsets of $G$, such that $G = G_1 \cup G_2$. Write $A = (cl(x) \cup G_1) \cup (cl(x) \cup G_2)$. It is clear that $A$ is $\beta$-reducible. ∎

The following result is immediate consequence of Lemma 4.6.

**Corollary 4.7** *Let $\tau \in \mathcal{L}_\rho$ be presober. Then every $\beta$-irreducible subset of $X$ is also $\tau$-irreducible.*

**Proposition 4.8** *Let $\tau \in \mathcal{L}_\rho$ be presober, $x \in X$ and $\mathcal{F}$ a filter on $X$. If $\mathcal{F} \supseteq N_{at(\rho)}(x)$ then $\beta = \tau \cap (\mathcal{E}(x) \cup \mathcal{F})$ is presober.*

**Proof.** If $\mathcal{F} \supseteq N_{at(\rho)}(x)$, then $\beta \in \mathcal{L}_\rho$ (Lemma 3.5(i)). Given a $\beta$-irreducible set $A$, then $A$ is $\tau$-irreducible (Corollary 4.7), and hence $A$ is the $\tau$-closure of a point, and thus the $\beta$-closure of a point. Therefore $\beta$ is presober. ∎

**Proposition 4.9** *There are not minimal presober members of $\mathcal{L}_\rho$.*

**Proof.** Let $\tau \in \mathcal{L}_\rho$ be a presober topology. Since $at(\rho)$ can not be presober, there is $V \in \tau \backslash at(\rho)$. Let $y \in V$ and let $\beta = \tau \cap (\mathcal{E}(y) \cup N_{at(\rho)}(y))$. By Proposition 4.8, $\beta$ is a presober topology which is obviously strictly weaker than $\tau$. Therefore $\tau$ is not minimal presober. ∎

**Corollary 4.10** *There are not minimal (sober and $T_1$) topologies on an infinite set.*

**Proof.** Follows from Proposition 4.9 with $\rho$ any $T_1$-topology. ∎

# References

[1] O.T. Alas, R.G. Wilson, *Minimal properties between $T_1$ and $T_2$*, Houston J. Math., 32(2)(2006), 493-504.

[2] A. Bella, C. Constantini, *Minimal KC spaces are compact,* Topology Appl., 155 (13)(2008), 1426-1429.

[3] M.P. Berri, M*inimal topological spaces*, Trans. Amer. Mat. Soc., Vol 108 (1)(1963), 97-105.

[4] M.P. Berri, R.H. Sorgenfrey, *Minimal regular spaces*, Proc. Amer. Math. Soc., Vol. 14 (3)(1963), 454-458.

[5] M.P. Berri, Hon-Fei Lai, D.H. Pettey, *Non compact, minimal regular spaces*, Pacific J. Math., Vol. 69 (1)(1977), 19-23.

[6] M.L. Colasante, C. Uzcátegui, J. Vielma, *Boolean algebras and low separation axioms.* Top. Proc., Vol. 34(2009), 1-15.

[7] A.S. Davis, *Indexed systems and neighborhoods for general topological spaces*, Amer. Math. Monthly, (68)(1961), 886-893.

[8] J. Dontchev, M. Ganster, G.J. Kennedy, S.D. McCartan, *On minimal Door, minimal anti-compact and minimal $T_{3/4}$-spaces*, Math. Proc. Royal Irish Acad.. Vol. 98A (2)(1998), 209-215.

[9] J.E. Joseph, *Characterizations of minimal Hausdorff Spaces*, Proc. Amer. Math. Soc., Vol. 61 (1)(1976), 145-148.

[10] R.L. Larson, *Minimal $T_0$-spaces and Minimal $T_D$-spaces*, Pacific J. Math., Vol. 31 (2)(1969), 451-458.

[11] A.E. McCluskey, S.D. McCartan, Minmality with respect to Young's axiom, Houston J. Math. 21(2)(1995), 413-428.

[12] J. Porter, J. Thomas, *On H-closed minimal Hausdorff spaces*, Trans. Amer. Math. Soc., Vol. 138 (1969), 159-170.

[13] C.T. Scarborough, Minimal Urysohn spaces, Pacific J. Math., 3(1968), 611-617.

[14] A.N. Shanin, *On separation in topological spaces*, Doklady SSSR, (38)(1943), 110-113.

[15] M.K. Singal, Asha Mathur, *A note on minimal spaces*, Israel J. Math., Vol. 12(2)(1972), 204-206.

[16] J. Stephräns, S. Watson, *Mutually complementary families of $T_1$ topologies, equivalence relations and partial orders*, Proc. Amer. Math. Soc., Vol. 123 (7)(1995), 2237-2249.

[17] C. Uzcátegui, J. Vielma, *Alexandroff topologies viewed as closed sets in the Cantor cube*, Divulg. Mat., (13)(2005), 45-53.

[18] G. Viglino, *A co-topological application to minimal spaces*, Pacific J. Math., Vol. 27(1)(1968), 197-200.

[19] S. Willard, *General Topology*, Addison-Wesley Publishing Co., 1970.

**MARÍA L. COLASANTE**
Departamento de Matemáticas, Facultad de Ciencias,
Universidad de Los Andes
Mérida 5101, Venezuela
e-mail: marucola@ula.ve

**D. Van der ZYPEN**
Allianz Suisse Insurance, Bleicherweg 19
CH-8022 Zurich, Switzerland
e-mail: dominic.zypen@gmail.com